\documentclass[11pt, a4paper]{article}
\usepackage{amsmath}
\usepackage{amsthm}
\usepackage{amsfonts}
\usepackage{amssymb,amsmath,latexsym}
\usepackage{graphicx}
\input epsf

\usepackage{epsfig}    
\usepackage{bbm}       

\usepackage[latin1]{inputenc}

\usepackage{epic}      
\usepackage{eepic}     
\usepackage{calc}      

\usepackage{xy}        
\xyoption{all}
\theoremstyle{plain}
\theoremstyle{definition}

\newtheorem{aussage}{Aussage}[section]

\newtheorem{thm}[aussage]{Theorem}

\newtheorem{prop}[aussage]{Proposition}
\newtheorem{cor}[aussage]{Corollary}
\newtheorem{rem}[aussage]{Remark}
\newtheorem{defn}[aussage]{Definition}

\newtheorem{claim}[aussage]{Claim}

\numberwithin{equation}{section}

\newcommand{\cyc}{\operatorname{cyc}}

\newcommand{\res}{\operatorname{res}}

\newcommand{\mb}[1]{\mathbb{#1}}
\newcommand{\mf}[1]{\mathfrak{#1}}

\newcommand{\bC}{\mathbb{C}}

\newcommand{\bE}{\mathbb{E}}

\newcommand{\bN}{\mathbb{N}}

\newcommand{\bR}{\mathbb{R}}

\newcommand{\bZ}{\mathbb{Z}}

\newcommand{\cA}{\mathcal{A}}
\newcommand{\cB}{\mathcal{B}}

\newcommand{\cD}{\mathcal{D}}
\newcommand{\cE}{\mathcal{E}}

\newcommand{\cG}{\mathcal{G}}
\newcommand{\cH}{\mathcal{H}}

\newcommand{\cK}{\mathcal{K}}
\newcommand{\cL}{\mathcal{L}}
\newcommand{\cM}{\mathcal{M}}
\newcommand{\cN}{\mathcal{N}}
\newcommand{\cO}{\mathcal{O}}
\newcommand{\cP}{\mathcal{P}}
\newcommand{\cQ}{\mathcal{Q}}
\newcommand{\cR}{\mathcal{R}}
\newcommand{\cS}{\mathcal{S}}
\newcommand{\cT}{\mathcal{T}}
\newcommand{\cU}{\mathcal{U}}

\newcommand{\cX}{\mathcal{X}}

\begin{document}

\title{The Hilbert-uniformization is real-analytic}
\date{}
\maketitle
\begin{center}
{Johannes F. Ebert$\footnote{e-mail: ebert@math.uni-bonn.de}$ and   
Roland~M. Friedrich$\footnote{ e-mail: rolandf@mpim-bonn.mpg.de}$ 
} \\

\vspace*{3mm}
{\em Max-Planck-Institut f\"ur Mathematik, Vivatsgasse 7, D-53111 Bonn} \\
\vspace*{6mm}

\end{center}

\begin{abstract}
In \cite{Boed}, C.-F. B\"odigheimer constructed a finite cell-complex $\mf{Par}_{g,n,m}$ and a bijective map $\cH: \mf{Dip}_{g,n,m} \to \mf{Par}_{g,n,m}$ (the Hilbert-uniformization) from the moduli space of dipole functions on Riemann surfaces with $n$ directions and $m$ punctures to $\mf{Par}_{g,n,m}$. In \cite{Boed} and \cite{Eb}, it is proven that $\cH$ is a homeomorphism.
The first result of this note is that the space $\mf{Dip}_{g,n,m}$ carries a natural structure of a real-analytic manifold.
Our second result is that $\cH$ is real-analytic, at least on the preimage of the top-dimensional open cells of $\mf{Par}_{g,n,m}$.
\end{abstract}

\section{Introduction}

\label{}
Throughout this note, $F$ will denote a closed smooth oriented surface of genus $g$. The genus can be arbitrary, but we will present the arguments in detail only in the case $g \geq 2$. The case $g=0,1$ can be treated with minor changes. Furthermore, let $\cP,\cQ \subset F$ be finite sets of points, $\cP \cap \cQ= \emptyset, \cP=\{P_1, \ldots , P_m\}; \cQ=\{Q_1, \ldots , Q_n\}, \cX=\{X_{1}, \ldots ,X_{n}\}$, $X_j \in T_{Q_j}F \setminus \{0\}$ are nonzero tangent vectors. Assume that $P_i \neq P_j$ if $i \neq j$ and the same for the points $Q_i$. We will always assume that $n>0$, i.e., we have at least one marked point with a tangent vector. 

\begin{defn} Let $(F,\cP,\cX)$ be as above and let a complex structure $J$ on $F$ be given. A \emph{dipole function} with singularities $(\cX,\cP)$ is a harmonic function $u: F \setminus \cP \cup \cQ \to \mathbb{R}$ with singularities as follows:
\begin{enumerate}
\item If $P_i \in \cP$ and $z$ is a local holomorphic coordinate with $z(P_i)=0$, then $u(z)=C_i \log(|z|)+b(z)$ where $C_i>0$ and $b$ harmonic.
\item If $Q_j \in \cQ$ and $z$ is a local holomorphic coordinate with $Tz(X_j)=\frac{\partial}{\partial x}\big|_{0} \in T_0 \mb{C}$, then $u(z)=\Re(\frac{1}{z}) + B_j \log(|z|) +a(z)$, where $a$ is harmonic.
\end{enumerate}
\end{defn}

Note that the numbers $C_i$ and $B_i$ do not depend on the choice of the coordinate and that their sum must be zero, due to the residue theorem. We will only consider the complex differential $\omega= \partial u$ of a dipole function. The possible $\omega $ are characterized by the conditions:
\begin{enumerate}
\item $\omega$ is a meromorphic 1-form on $F$, which has poles of order $1$ at the points of $P$ and poles of order $2$ at the points of $Q$, furthermore,
\item the periods of $\omega$ over all homology classes in $H_1(F \setminus P \coprod Q;\bZ)$ are purely imaginary;
\item the residues at the points of $P$ are positive real numbers, the residues at the points of $Q$ are arbitrary real numbers.
\end{enumerate}

\begin{defn}\label{definition1}
The space $\mathfrak{Dip}_{g,n,m}$ as a set is defined to be the set of equivalence classes of quadruples $(F,\cP,\cQ,\omega)$ , where two such objects $(F,\cP,\cQ,\omega)$ and $(F^{\prime},\cP^{\prime},\cQ^{\prime},\omega^{\prime})$ are considered to be equivalent if there exists a biholomorphic map $f:F \to F^{\prime}$ with $f(P_i)=P_{i}^{\prime}$ and $f^{*}(\omega^{\prime})=\omega$.
\end{defn}

\begin{rem}
The space $\mf{Dip}_{g,n,m}$ is almost the same as the space $\mf{Pot}^{\bC}(g,n,m)$, which is studied in detail in \cite{Eb}, \cite{Boed}. In fact, a differential $\omega$ as in Definition~\ref{definition1} can be written as $\partial u$, where $u$ is a real-valued harmonic function with singularities, a so-called dipole function. $u$ is unique up to addition of a constant. Hence: $\mf{Dip}_{g,n,m} \cong \mf{Pot}^{\bC}(g,n,m) / \bR$, the factor $\bR$ referring to addition of real constants.
\end{rem}

A construction of a topology on $\mathfrak{Dip}_{g,n,m}$ is given in \cite{Eb}, using the fiber bundle Teichm\"uller theory of Earle and Eells (\cite{EaEe}) and methods from global analysis. We will give a different construction of this topology, which also puts a real-analytic structure on $\mathfrak{Dip}_{g,n,m}$. This construction relies on more advanced results from Teichm\"uller theory and complex analysis. Our first main result is the theorem below.

\begin{thm}\label{theoremone}
The space $\mf{Dip}_{g,n,m}$ is a real-analytic manifold of (real) dimension $6g-6+5n+3m-1$.
\end{thm}

B\"odigheimer (\cite{boed2,Boed}) constructs a finite cell complex $\mf{Par}_{g,n,m}$ and two mutually inverse maps $\cH:\mf{Dip}_{g,n,m} \to \mf{Par}_{g,n,m}$ and $\cG:\mf{Par}_{g,n,m} \to \mf{Dip}_{g,n,m}$. The first author of the present paper showed that both maps are continuous. In this note, we will show that $\cH$ is even real-analytic.\\
Next, we will describe the open, top-dimensional cells of the complex $\mf{Par}_{g,n,m}$. The combinatorics becomes much more complicated for the lower-dimensional cells and is not needed in this paper, since we only consider the maps $\cG$ and $\cH$ on these top-dimensional cells.\\
For simplicity, we assume that $n=1$ in the rest of the paper. Set $h:=2g+m$. A top dimensional cell is given by the following combinatorial data:

\begin{enumerate}
\item A sequence of permutations $\sigma_0, \ldots , \sigma_h \in \Sigma_{2h+1}$, such that $\sigma_i (2h) =0$ (the symmetric group $\Sigma_{2h+1}$ acts on the set $\{0, \ldots ,2h \}$) and such that $\sigma_0(j) = j+1$ for all $j$ except $j=2h$.
\item A bijection $\nu: \{0,1, \ldots , m\} \to \cyc(\sigma_h)$, where $\cyc(\sigma)$ denotes the set of all cycles of the permutation $\sigma$, such that $\nu(0) $ is the cycle containing $0$.
\end{enumerate}

This cell is $h-1 + 2h-1$-dimensional (which is the same as $\dim \mf{Dip}_{g,1,m}$). The coordinates of a point in the cell are given by two sequences of real numbers

\begin{enumerate}
\item $a_h < a_{h-1}< \ldots < a_1$,
\item $b_1 < b_2 < \ldots < b_{2h} $,
\end{enumerate}

where two sequences are identified if the sequence of $a$'s and $b$'s differs by real constants. We call a point in the complex $\mf{Par}_{g,1,m}$ a \emph{ parallel slit domain} of type $(g,1,m)$.\\
Given such data, we can construct a Riemann surface together with a dipole function using the following procedure. Set $a_{h+1}= b_0 = - \infty $ and $a_{0}= b_{2h+1}= + \infty $ and define rectangles $R_{i,j} := [a_{i+1},a_i] \times [b_{j-1}, b_j] \subset \bR^2 = \bC$. Now we identify the left boundary of the rectangle $R_{i,j}$ with the right boundary of $R_{i-1,j}$ (for $i=0, \ldots , h$, all $j$) and the upper boundary of $R_{i,j}$ with the lower boundary $R_{i,\sigma_i (j)}$ ($j=0, \ldots ,h$, all $i$). The result is a topological surface of genus $g$ with $m+1 $ punctures which has a canonical complex structure. It turns out that the punctures can be filled in to obtain a closed Riemann surface. The function $z \mapsto \Re (z)$ clearly gives a dipole function on the resulting closed surface. The whole procedure yields the map $\cG$.\\
For the construction of the inverse map, we need the next

\begin{defn}
Let $\omega=\partial u$, where $u$ is a real-valued dipole function. Then the \emph{critical graph} $\cK$ of $\omega$ is the union of all flow lines of the ascending gradient flow of $u$ whose limit for $t \to + \infty$ is a critical point.\\
$\omega$ is called \emph{generic} if all zeroes are simple, if the values of $u$ at the critical points are all distinct and if there are no integral curves connecting two critical points.
\end{defn}

The complement of the critical graph is contractible (if $n=1$, otherwise it is the union of $n$ contractible components), thus the function $u$ is the real part of a holomorphic function $f$ (or, equivalently, $df= \partial u$).
It turns out that $f$ is injective. From the image of the function $f$, one can construct a parallel slit domain. This gives the map $\cH$, which is also called the \emph{Hilbert-uniformization}. We describe how the numbers $a_i$ and $b_j$ are constructed. Let $x_1, \ldots ,x_h$ be the critical points of $u$. For any pair of critical points $x_k, x_l$, there is a unique (up to homotopy) path $c_{k,l} $ in $F$ which does not meet $\cK$ except at the endpoints. By integrating $\omega$ over $c_{k,l}$, we obtain a complex number $z_{k,l}$. Note that $z_{k,k}$ is purely imaginary (and nonzero). Choose $a_1$ and $b_1$ arbitrarily. By adding the different numbers $z_{k,l}$, we obtain $2h$ complex numbers, which are arranged in pairs of numbers with coinciding real part.
Separate the real and imaginary parts and order the numbers as required in the definition of a parallel slit domain.\\
It turns out that generic dipole functions correspond to points in the top-dimensional cells.\\
Now we state the second main result of the present paper.

\begin{thm}\label{theoremtwo}
Let $\cR \subseteq \mathfrak{Dip}_{g,m,n}$ be the open subset consisting of generic differentials\footnote{The openness of this map follows from the results of \cite {Boed} and \cite{Eb}}. Then the restriction $\mathcal{H}|_{\cR}$ is real-analytic.
\end{thm}

\section{The moduli space of dipole functions}\label{modulr}

In this section, we prove Theorem~\ref{theoremone}. We start with a general folklore result from complex algebraic geometry.

\begin{prop}\label{family}
Let $E , B$ be complex manifolds and $\pi: E \to B$ be a holomorphic proper submersion (hence a differentiable fiber bundle with compact complex manifolds as fibers, see \cite{Voi}). Let $L \to B$ be a holomorphic line bundle. Assume that the function $B \ni b \mapsto \dim H^0(E_b, \cO (L_b)) \in \bN$ is constant ($E_b:= \pi^{-1}(b)$; $L_b:= L|_{E_b}$.). Then the disjoint union $V:=\bigcup_{b \in B} H^0(E_b; \cO (L_b))$ has a unique structure of a holomorphic vector bundle on $B$, such that the evaluation map $V \times_{B} E \to L$ is holomorphic.
\end{prop}

From chapter 5 in the book \cite{Nag}, one extracts the following facts: Let $\cT:=\cT_{g,n+m}$ be the Teichm\"uller space of genus $g$ Riemann surfaces with $n+m$ distinct distinguished points. This is a complex manifold of dimension $3g-3+n+m$. Further there exists a complex manifold $\cE:=\cE_{g,n+m}$ and a fiber bundle $\pi: \cE_{g,n+m} \to \cT_{g,n+m}$ and $m+n$ holomorphic cross-sections $s_1, \ldots ,s_{m+n}$. The maps $\pi,s_i$ are holomorphic. Let $x:=[(F,p_1, \ldots ,p_{n+m})] \in \cT_{g,n+m}$. Then $\pi^{-1}(x)$ is a Riemann surface, canonically isomorphic to $F$; and the points $s_i(x)$ correspond to the points $p_i$. In other words, $\pi$ is a tautological surface bundle. It is universal among homotopically trivialized families of Riemann surfaces with constant punctures.\\
The mapping class group $\Gamma_{g,n+m}$ acts on both $\cE_{g,n+m}$ and $\cT_{g,n+m}$; $\pi$ is an equivariant map and both actions are holomorphic\footnote{This follows from the statement in chapter 5.5 in \cite{Nag}} and properly discontinuous. The action is \emph{not} free, but it has finite stabilizers and, if $n+m >0$, the stabilizers are cyclic groups.\\
Let $\Omega_{\cE/\cT}$ be the holomorphic vertical cotangent bundle of $\cE$. Let $S_i$ be the image of $s_i$. It is a complex submanifold of $\cE$. Define a divisor on $\cE$ by $\cD:= 2S_1 + \ldots +2S_n + S_{n+1} + \ldots + S_{n+m}$. $\cD$ defines a holomorphic line bundle $[\cD]$ on $  \cE$. Define $\cL:= \Omega_{\cE/\cT} \otimes [\cD]$. Note that $\cL$ has a $\Gamma_{g,n+m}$-action which is holomorphic. For any point $x \in \cT$, let $\cL_x$ be the restriction of $\cL$ to the fiber $\cE_x:=\pi^{-1}(x)$ over $x$. It is a holomorphic line bundle over the Riemann surface representing $x$ and its global holomorphic sections are the meromorphic 1-forms satisfying the condition (iv) in Definition~\ref{definition1}, where $p_i = s_i(x)$.
Let $\cA$ be the disjoint union $\bigcup_{x \in \cT} H^0(V_x; \cO (L_x))$. It is a holomorphic vector bundle on $\cT$ by Proposition \ref{family}. The complex dimension of its total space is $3g-3+3n+2m-1$. $\cA$ comes with a holomorphic action of $\Gamma_{g,n+m}$, which turns it into a $\Gamma_{g,n+m}$-equivariant vector bundle.\\

Now remember that $\cT_{g,n+m}$ is contractible; thus the bundle $\cE \to \cT$ can be topologically trivialized, in other words, there exists a diffeomorphism $\Psi:\cE \to \cT \times F$, commuting with the maps to $\cT$. Note that $\Psi$ cannot chosen to be holomorphic. The trivialization can be chosen to respect all the sections $s_1, \ldots ,s_{m+n}$. Let $c \in H_1(F \setminus \{p_1, \ldots ,p_{n+m}\})$ be a homology class and choose a path in $F$ representing $c$. Also, choose paths around the points $p_i$. 

\begin{cor}\label{residueperiod}
\begin{enumerate}
\item The map $\res:\cA \to \cT \times \bC^{m+n-1}$,\\ $\omega \mapsto (\rho(\omega),((\res_{p_1}\omega, \ldots ,\res_{p_m}\omega, \res_{q_1}\omega ,\ldots ,\res_{q_{n-1}} \omega))$ is a holomorphic vector bundle epimorphism.
\item For any homology class $c \in H_1(F \setminus P \cup Q, \bZ)$, the period map $\cA \to \cT \times \bC$; $\omega \mapsto (\rho(\omega),\int_{c} \omega)$, is holomorphic.
\end{enumerate}
\end{cor}

\textbf{Proof:} This follows immediately from Proposition \ref{family} and from standard results of complex analysis, because the residuum and the period maps are given by the evaluation map and by integration.
\qed

\begin{defn}
Let $\cB \subset \cA$ be the subspace of all meromorphic forms $\omega$, satisfying:
\begin{enumerate}
\item all periods of $\omega$ are purely imaginary,
\item all residues at the points $p_1, \ldots, p_{n+m}$ are real,
\item the residues at $p_{n+1}, \ldots ,p_{m+n}$ are positive,
\item at the points $p_1, \ldots, p_n$, $\omega$ has a pole of order two.
\end{enumerate}
Let $\cB_0 \subset \cB$ be the subset consisting of the differentials whose zeroes are all simple and whose critical graphs are generic (see section \ref{hilbertuniform} for a definition). Let also $\cA_0$ be the space of meromorphic forms in $\cA$ such that all zeroes are simple. Then $\cB_0 = \cB \cap \cA_0 $.
\end{defn}

The conditions are real-analytic conditions, such that $\cB$ and $\cB_0$ are real-analytic submanifolds of $\cA$. Further, the defining conditions are $\Gamma_{g,n+m}$-equivariant. The quotient $\cB / \Gamma_{g,n+m}$ is the space $\mathfrak{Dip}_{g,n,m}$ from Definition~\ref{definition1}. It is also a real-analytic manifold.\\
We also note that $\cB_0 / \Gamma_{g,n+m} = \cR$.\\
Using the results of \cite{EaEe}, it can be shown that the two ways to define the topology on $\mf{Dip}$ lead to the same topology.

\section{Hilbert-uniformization}\label{hilbertuniform}

The proof of Theorem~\ref{theoremtwo} is given in this section.

There is a diagram

$$\xymatrix{
 & \cN \ar[d] &        & \cD \ar[d] \\
  & \cU_0 \ar[r] \ar[d]^{\rho} \ar[ld]^{ev} & \cU \ar[r] \ar[d] & \cE \ar[d]^{\pi}\\
  T^{*}(\cU_0 / \cA_0) \ar[r] & \cA_0 \ar[r] & \cA  \ar[r]^{q} & \cT.\\
}$$

Here $\cU:=\cA \times_{\cT} \cE$ is the total space of the pullback of the surface bundle $\cE \to \cT$. $\cD$ is the divisor from section \ref{modulr}. $\cU_0$ is the complement of $\cD$ in the restriction of $\cU$ to the subspace $\cA_0$.\\
Recall that $\cA_0$ consists of those differentials whose zeroes are simple. Thus the evaluation map $ev$ is transverse to the zero section $\cS$ in the bundle $T(\cU_0 / \cA_0) \to \cU_0$.\\ 
By Proposition \ref{family} the map $ev$ is holomorphic. Hence:\\

\begin{cor}
$\cN := ev^{-1}(\cS)$ is a complex-analytic submanifold of complex codimension 1. Further, $\rho|_{\cN}\to \cA_0$ is a finite holomorphic covering.
\end{cor}

To finish the proof of Theorem~\ref{theoremtwo}, take an open subset $\cO \subset \cA_0$, where the covering $\cN \to \cA_0$ is trivial (this means, the zeroes of the differential can be numbered). For two distinct sheets $x_0, x_1$ of the covering $\cN$, take a homotopy class (relative $x_i$) of paths $c$ going from $x_0$ to $x_1$ (remember that the bundle $\cU_0 \to \cA_0$ is topologically trivial).\\
Then we obtain a function $\alpha_c:\cO \to \bC$ by the rule $\omega \mapsto \int_{c} \omega$. $\alpha_c$ is holomorphic by a reasoning analogous to the proof of corollary \ref{residueperiod}.\\
Now we restrict our attention to the real-analytic submanifold $\cB_0$ of $\cA_0$. Let $\cM:= \cN \cap \rho^{-1}(\cB_0)$.\\
$\cM$ consists of all critical points of genuine dipole functions. The homotopy class $c$ should be chosen, such that the path does not meet the critical graphs except at the endpoints (compare the definition of the Hilbert-uniformization). The function $\alpha_c$ gives the coordinates of the Hilbert-uniformization. Since it is the restriction of a holomorphic function to a real-analytic submanifold, it is real-analytic.
\qed

\section{Concluding remarks} 

It seems reasonable that the inverse map $\cG$ is also real-analytic, at least on the top-dimensional cells. One way to prove this is to compute the differential of $\cH$ and to prove that it is nonsingular. But due to the lack of a suitable coordinate system on $\mf{Dip}_{g,n,m}$, we postpone this.\\
However, there is a partial result, which can be obtained quite easily. First of all, the residues of the resulting dipole function obviously depend real-analytically on the coordinates of the parallel-slit domain.
Secondly, consider the projection $\cB \to \cT$. There is a criterion to determine whether a map into Teichm\"uller space is real-analytic, namely:

\begin{claim}\label{ahlbers} Let $B$ be a real-analytic manifold and let $p:E \to B$ be a topologically trivialized family of Riemann surfaces of genus $g$. Let $F_0$ be a closed Riemann surface and let $\Phi: B \times F_0 \to E$ a continuous map commuting with the projections to $B$. Assume that for all $b \in B$, the restriction $\Phi|_{ \{b\} \times F_0} : \{b\} \times F_0 \to p^{-1}(b)$ is a quasiconformal homeomorphism of Riemann surfaces. Assume that the complex dilatation $\mu:B \to L^{\infty}( F_0 \to \bE)$ (set $\bE:=\{|z| <1 \}$ is a real-analytic function. Then the moduli map of the family $B  \to \cT_g$ is real-analytic (with respect to the natural real-analytic structure on $\cT_g$).
\end{claim}

This should follow quite easily from the results of Ahlfors and Bers on Teichm\"uller theory. One can apply this criterion in the given situation and we obtain that the composition $\mf{Par}_{g,n,m}\supset C \to \cB \to \cT$ is real-analytic on a top-dimensional cell $C$ (here we use a lift of $C\to \mf{Dip}_{g,n,m}$ to the cover $\cB$ of $\mf{Dip}_{g,n,m}$).
The remaining problem is to find a coordinate which detects the direction of the tangent vector. Because the Teichm\"uller maps are not differentiable, classical Teichm\"uller theory does not seem to be the appropriate framework for an approach to this problem.

\section*{Acknowledgements}

The authors would like to acknowledge the hospitality and the support of the Max-Planck-Insitute for Mathematics in Bonn. We thank C.-F. B\"odigheimer for discussions. Also we thank D. Huybrechts for his remark concerning Prop.~\ref{family}.

\end{document}